\documentclass[11pt]{article}

\usepackage{latexsym,amssymb}

\usepackage{amscd}

\input xypic
\xyoption{all}

\newcommand{\rar}{\rightarrow}
\newcommand{\lar}{\longrightarrow}

\newtheorem{theorem}{Theorem}[section]

\newtheorem{corollary}[theorem]{Corollary}
\newtheorem{proposition}[theorem]{Proposition}

\newtheorem{example}[theorem]{Example}
\newtheorem{conjecture}[theorem]{Conjecture}

\def\demo{\noindent{\bf Proof. }}

\def\text{\mbox}

\newcommand{\Div}{\mbox{\rm Div}}
\newcommand{\Cl}{\mathsf{Cl}}
\newcommand{\Prin}{\mbox{\rm Prin}}

\newcommand{\Ext}{\mbox{\rm Ext}}
\newcommand{\Hom}{\mbox{\rm Hom}}

\newcommand{\ann}{\mbox{\rm ann}}
\newcommand{\gr}{\mbox{\rm gr}}

\newcommand{\Rees}{\mbox{${\mathcal R}$}}
\def\QED{\hfill$\Box$}
\def\spec{\mbox{\rm Spec}}
\def\ddet0{\mbox{\rm det}_0 }
\def\div{\mbox{\rm div} }
\newcommand{\coker}{\mbox{\rm coker }}
\newcommand{\Gr}{\mathrm G}
\def\gr{\mbox{\rm gr}}

\newcommand{\emb}{\mbox{\rm edim}}

\newcommand{\rank}{\mbox{\rm rank }}
\newcommand{\codim}{\mbox{\rm codim}}

\newcommand{\kernel}{\mbox{\rm ker}}
\newcommand{\im}{\mbox{\rm image}}

\def\depth{\mbox{\rm depth }}
\def\grade{\mbox{\rm grade }}
\def\height{\mbox{\rm height}}

\begin{document}

\title{\LARGE\bf{\sc
Rees Algebras of Conormal Modules} }

\author{
{\normalsize\sc Jooyoun Hong}\\
{\small Department of Mathematics, Rutgers University, Piscataway, NJ 08854-8019, USA}\\
{\small e-mail: jooyoun@math.rutgers.edu}}

\maketitle

\begin{abstract}
We deal with classes of prime ideals whose associated graded ring
is isomorphic to the Rees algebra of the conormal module in order
to describe the divisor class group of the Rees algebra and to
examine the normality of the conormal module.
\end{abstract}


\footnotetext{Correspondence to: Jooyoun Hong, Department of
Mathematics, Purdue University, West Lafayette, IN 47907-2067,
USA. e-mail: jooyoun@math.rutgers.edu
\\ Mathematics Subject Classification 2000: 13H10, 13B21, 13C20.
\\ This paper is based on the Ph.D. dissertation of the author at Rutgers University, written under the direction of
Professor Wolmer V. Vasconcelos.}

\section{Introduction}\label{Intro}

We investigate several algebras associated to a prime ideal of a
commutative Noetherian ring and study the relationships among them
by examining normality conditions (divisor class groups,
the construction of integral closures, etc). Throughout this
paper, let $R$ be a commutative Noetherian ring and $\mathfrak{p}$
a prime ideal of $R$. We illustrate the algebras to be examined
with the following diagram (unexplained terminology will be
discussed in the text):

\bigskip

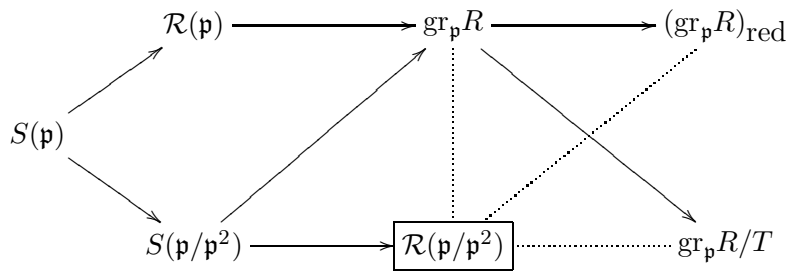
\[
\diagram
   &\Rees(\mathfrak{p}) \rrto  && \gr_{\mathfrak{p}}
R \dddotted \rrto \ddrrto &&(\gr_{\mathfrak{p}} R)_{\mbox {\rm
red}}
 \ddlldotted\\
S(\mathfrak{p}) \urto \drto  &&& && \\
   &S(\mathfrak{p}/\mathfrak{p}^2)  \rrto \uurrto
   &&\Rees(\mathfrak{p}/\mathfrak{p}^2)\framed
   && \gr_{\mathfrak{p}} R/T \lldotted \\
\enddiagram
\]

\noindent where $S(\mathfrak{p})$ is the symmetric algebra of
$\mathfrak{p}$, $S(\mathfrak{p}/\mathfrak{p}^2)$ is the symmetric
algebra of $\mathfrak{p}/\mathfrak{p}^2$ as an
$R/\mathfrak{p}$--module, $\Rees(\mathfrak{p})$ is the Rees
algebra of $\mathfrak{p}$, $\gr_{\mathfrak{p}} R$ is the
associated graded ring, $(\gr_{\mathfrak{p}} R)_{\mbox{ \rm red}}$
is the reduced ring of the associated graded ring,
$(\gr_{\mathfrak{p}} R)/T$ is the associated graded ring mod its
torsion $T$, and $\Rees(\mathfrak{p}/\mathfrak{p}^2)$ is the Rees
algebra of $\mathfrak{p}/\mathfrak{p}^2$, which is the focus of
our study.

Let us recall the notion of the Rees algebra of a module. Let $R$
be a domain and $E$ a finitely generated torsionfree $R$--module
of rank $e$ with an embedding $\varphi: E \hookrightarrow R^e$.
The {\em Rees algebra} $\Rees(E)$ of $E$ is the subalgebra of the
polynomial ring $R[T_1, \ldots, T_e]$ generated by all linear
forms $a_1T_1 + \cdots + a_eT_e$, where $(a_1 , \ldots , a_e)$ is
the image of an element of $E$ in $R^e$ under the embedding. The
Rees algebra $\Rees(E)$ is a standard graded algebra
$\bigoplus_{n=0}^{\infty} E_n$ over $R$ with $E_1 =E$. We refer
the reader to \cite{EHU1} for a general discussion of Rees
algebras of modules over general rings, including the fact that
the Rees algebra $\Rees(E)$ is independent of the embedding
$\varphi$ when $R$ is a domain. In general, there is a surjection
from the symmetric algebra of $E$ onto the Rees algebra of $E$ and
the module $E$ is said to be of {\em linear type} if this
surjection is an isomorphism. Let $U \hookrightarrow R^e$ be a
submodule of $E$ with the same rank as that of $E$. The module $E$
is {\em integral} over the module $U$ if the Rees algebra of $E$
is integral over the $R$-subalgebra generated by $U$. In this case
we say that $U$ is a {\em reduction} of $E$. The {\em integral
closure} $\overline{E}$ of $E$ is the largest submodule of $R^e$
which is integral over the module $E$. If $E$ is equal to
$\overline{E}$, then $E$ is called {\em integrally closed} or {\em
complete}. If the Rees algebra $\Rees(E)$ of $E$ is integrally
closed, then the module $E$ is said to be {\em normal}. This means
that each component $E_n$ of $\Rees(E)$ is integrally closed
(Proposition~\ref{intclosrees}).

We focus on the normality of a conormal module
$\mathfrak{p}/\mathfrak{p}^2$ as an $R/\mathfrak{p}$--module. In
general, the Rees algebra $\Rees(\mathfrak{p}/\mathfrak{p}^2)$
depends only on the module $\mathfrak{p}/\mathfrak{p}^2$ over the
ring $R/\mathfrak{p}$, not on the ring $R$ itself. To ensure $R$
has a role, we must force a relationship between $\mathfrak{p}$
and $R$, which is the case when we assume, for example, that the
prime ideal $\mathfrak{p}$ has finite projective dimension over
$R$. Throughout this paper, we denote the associated graded ring
of a prime ideal $\mathfrak{p}$ by $\Gr$, i.e.,
\[\Gr = \gr_{\mathfrak{p}} R = \bigoplus_{t=0}^{\infty} \mathfrak{p}^t /\mathfrak{p}^{t+1}=\bigoplus_{t=0}^{\infty} \Gr_t.\]

Our goal is to describe the divisor class group of an integrally
closed associated graded ring and to examine the normality of the
conormal module whose Rees algebra is isomorphic to the associated
graded ring. First, when the associated graded ring is integrally
closed, we obtain the following theorem.

\bigskip

\noindent {\bf Theorem~\ref{hong1} } {\it Let $R$ be a
Cohen--Macaulay ring and $\mathfrak{p}$ a prime ideal of finite
projective dimension. If $R/\mathfrak{p}$ is an integrally closed
domain of dimension greater than or equal to $2$ and the
associated graded ring $\Gr$ of $\mathfrak{p}$ is integrally
closed, then the mapping of divisor class groups
\[ \Cl(R/\mathfrak{p})\ni [L] \mapsto  [L\Gr] \in \Cl(\Gr)\]
is a group isomorphism. In particular, if $\Gr$ is integrally
closed and $R/\mathfrak{p}$ is a factorial domain, then  $\Gr$ is
a factorial domain.}

\bigskip

The exact sequence of the components of the embedding $\Gr \subset
\overline{\Gr}$,
\[ 0 \rar \Gr_t \lar \overline{\Gr}_t \lar C_t \rar 0, \quad t\geq 1,\]
raises the question of when the normality of the associated graded
ring $\Gr$ can be detected in the vanishing of $C_t$ for low $t$.
Thus, we study the relationship between the normality of the
conormal module and the completeness of components of the
associated graded ring. So far we prove the following.

\bigskip

\noindent {\bf Theorem~\ref{hong3} } {\it Let $R$ be a Gorenstein
local ring and $\mathfrak{p}$ a prime ideal generated by a
strongly Cohen--Macaulay $d$--sequence. Suppose that
$\mathfrak{p}$ has finite projective dimension, that
$R/\mathfrak{p}$ is an integrally closed domain of dimension $2$
and that the associated graded ring $\Gr$ of $\mathfrak{p}$ is a
domain. Then the conormal module $\mathfrak{p}/{\mathfrak{p}}^2$
is integrally closed if and only if $\Gr$ is normal.}

\bigskip

This paper is organized as follows. In Section~\ref{RAGR}, we give
a suitable condition for an associated graded ring to be a domain
or an integrally closed domain. Our focus is on classes of ideals
whose associated approximation complexes are acyclic. In
Section~\ref{DCG}, we prove our main result regarding the
relationship between the divisor class group $\Cl(\Gr)$ of the
associated graded ring $\Gr$ and that $\Cl(R/\mathfrak{p})$ of
$R/\mathfrak{p}$. In Section~\ref{ICCM}, we deal with special
kinds of prime ideals whose associated graded ring $\Gr$ is
isomorphic to the Rees algebra of the conormal module associated
to the prime ideal. We also study the normality of $\Gr$ in terms
of the completeness of the components of $\Gr$. More precisely, we
would like obtain any statement of the form $\Gr_t =
\overline{\Gr_t}$ for $t \leq n_0$ implies completeness for all
$t$. We will also pay attention to the required degrees of the
generators of the integral closure $\overline{\Gr}$ of the
associated graded ring $\Gr$.

\section{Rees Algebras as Associated Graded Rings}\label{RAGR}

In general, the Rees algebra of the conormal module
$\mathfrak{p}/\mathfrak{p}^2$ and the associated graded ring of
the prime ideal $\mathfrak{p}$ are distinct. They may even have
different Krull dimensions. For example, if $(R,\mathfrak{m})$ is
a local ring then the Rees algebra
$\Rees(\mathfrak{m}/\mathfrak{m}^2)$ is a polynomial ring with
$\nu(\mathfrak{m})$ number of variables over the residue field of
$R$, where $\nu(\mathfrak{m})$ is the minimal number of generators
of $\mathfrak{m}$. On the other hand, in the diagram illustrated
at the beginning, we suggest that there may be a relationship
between the associated graded ring of $\mathfrak{p}$ and the Rees
algebra of $\mathfrak{p}/\mathfrak{p}^2$. We point out when this
is possible.

\begin{proposition}\label{gci}
Let $(R,\mathfrak{m})$ be a universally catenary Noetherian local
ring, $\mathfrak{p}$ a prime ideal of $R$, and $\Gr$ the
associated graded ring of $\mathfrak{p}$. There exists a
surjection from $\Gr$ onto the Rees algebra
$\Rees(\mathfrak{p}/\mathfrak{p}^2)$ of the conormal module
$\mathfrak{p}/\mathfrak{p}^2$ if and only if $\mathfrak{p}$ is
generically a complete intersection (i.e., $R_{\mathfrak{p}}$ is a
regular local ring).
\end{proposition}

\demo If we assume that $\mathfrak{p}$ is generically a complete
intersection, then we apply the dimension formula of \cite[Lemma
1.2.2]{V2} to the Rees algebra
$\Rees(\mathfrak{p}/\mathfrak{p}^2)$ to obtain the following.
\[
\dim\Rees(\mathfrak{p}/\mathfrak{p}^2) \leq \dim R/\mathfrak{p} +
\emb(R_{\mathfrak{p}})= \dim R/\mathfrak{p} + \height(
\mathfrak{p}) = \dim \Gr,\] where $\emb(R_{\mathfrak{p}})$ is the
embedding dimension of $R_{\mathfrak{p}}$. This proves that
$\Rees(\mathfrak{p}/\mathfrak{p}^2)$ is a homomorphic image of the
associated graded ring $\Gr$.

Now suppose that there is a surjective homomorphism from $\Gr$ to
$\Rees(\mathfrak{p}/\mathfrak{p}^2)$, then
\[\emb(R_{\mathfrak{p}})=\dim \Rees(\mathfrak{p}_{\mathfrak{p}}/\mathfrak{p}_{\mathfrak{p}}^2) \leq \dim R_{\mathfrak{p}},\]
and hence $R_{\mathfrak{p}}$ is a regular local ring. \QED

\bigskip

Suppose that a prime ideal $\mathfrak{p}$ is of linear type. This
is equivalent to say that the natural map between the symmetric
algebra $\mathcal{S}(\mathfrak{p}/\mathfrak{p}^2 )$ and the
associated graded ring $\Gr$ of $\mathfrak{p}$ is an isomorphism.
The approximation complex $\mathcal{M}(\mathfrak{p})$ comes in at
this point since zeroth Koszul homology module
$H_0(\mathcal{M}(\mathfrak{p}))$ is
$\mathcal{S}(\mathfrak{p}/\mathfrak{p}^2)$. One condition on the
ideal $\mathfrak{p}$ that has an impact on the acyclicity of
$\mathcal{M}(\mathfrak{p})$ is called $G_{\infty}$: For each prime
ideal $\mathfrak{q}\in V(\mathfrak{p})$, the minimal number of
generators $\nu(\mathfrak{p}_{\mathfrak{q}})$ of
$\mathfrak{p}_{\mathfrak{q}}$ is less than or equal to height of
$\mathfrak{q}$. We are going to give a formulation of some
conditions leading to associated graded rings which are domains.
This formulation involves approximation complexes.

\begin{proposition}\label{appcx1}
Let $R$ be a Cohen--Macaulay local ring and $\mathfrak{p}$ a
non-maximal prime ideal with sliding depth. Then $\mathfrak{p}$
satisfies $G_{\infty}$ and the associated graded ring
$\Gr=\bigoplus G_t$ of $\mathfrak{p}$ is a domain if and only if
for every proper prime ideal $\mathfrak{q}$ containing $
\mathfrak{p}$,
\[\nu(\mathfrak{p}_{\mathfrak{q}})\leq \height(\mathfrak{q})-1,\]
where $\nu(\mathfrak{p}_{\mathfrak{q}})$ is the minimal number of
generators of $\mathfrak{p}_{\mathfrak{q}}$.
\end{proposition}

\demo Suppose that $\nu(\mathfrak{p}_{\mathfrak{q}})$ is less than
or equal to $\height(\mathfrak{q}) -1$ for every proper prime
ideal $\mathfrak{q}$ containing $\mathfrak{p}$. Then
$\mathfrak{p}$ satisfies $G_{\infty}$ and the approximation
complex $\mathcal{M}(\mathfrak{p})$ is acyclic (\cite[Theorem
5.1]{HSV1}). Let $S=\bigoplus S_t$ be the polynomial ring $R[T_1 ,
\ldots , T_n ]$ and $H_i$ the $i$th Koszul homology module
associated to $\mathfrak{p}$. Applying the Depth Lemma(\cite[Lemma
3.1.4]{V2}) to the approximation complex
\[ 0 \rightarrow H_r \otimes S_{t-r} \rightarrow \cdots \rightarrow H_0 \otimes S_t \rightarrow S_t(\mathfrak{p}/\mathfrak{p}^2)=G_t \rightarrow 0,  \]
gives that each component $\Gr_t$ of $\Gr$ satisfies
 \[\depth(\Gr_t) > \depth( H_r )-(r+1) \geq d-n+r -(r+1) \geq 0.\]
This implies that locally the algebra $\Gr$ is torsionfree as an
$R/\mathfrak{p}$--module. If $K$ is the field of fractions of
$R/\mathfrak{p}$, we have an embedding $\Gr \hookrightarrow
\Gr\otimes K$. Since $\Gr\otimes K$ is a ring of polynomials over
$K$, the associated graded ring $\Gr$ is a domain.

For the converse, we may assume that $(R,\mathfrak{q})$ is a local
ring. Since the ideal $\mathfrak{p}$ satisfies $G_{\infty}$, the
ideal $\mathfrak{p}$ is generated by a $d$--sequence and the
minimal number of generators $\nu(\mathfrak{p})$ is equal to the
analytic spread $\ell(\mathfrak{p})$ (\cite[Theorem 2.2]{Hu6}).
Therefore,

\[ \nu(\mathfrak{p}) = \ell(\mathfrak{p}) = \dim \Gr/\mathfrak{q}\Gr < \dim \Gr = \height(\mathfrak{q}),\]
where the inequality follows from the fact that $\Gr$ is a domain.
\QED

\begin{proposition} \label{appcx2} Let $R$ be a Cohen--Macaulay local ring and $\mathfrak{p}$ a prime ideal with sliding depth.
Suppose that $R/\mathfrak{p}$ is integrally closed. Then the
associated graded ring $\Gr$ of $\mathfrak{p}$ is an integrally
closed domain if for every proper prime ideal $\mathfrak{q}$
containing $\mathfrak{p}$,
\[\nu(\mathfrak{p}_{\mathfrak{q}})\leq \max\{ \height(\mathfrak{p}_{\mathfrak{q}}) \;,\; \height(\mathfrak{q})-2\}.
\] The converse holds if $\mathfrak{p}$ has finite projective
dimension.
\end{proposition}

\demo Let $\mathfrak{q'}$ be a prime ideal of $\Gr$ of height one
and $\mathfrak{q}$ its inverse image in $R$. The associated graded
ring $\Gr$ is Cohen-Macaulay (\cite[Theorem 5.1]{HSV1}) and
\[\Gr=\bigcap \Gr_{\mathfrak{q'}},\;\;\;\; \mbox{\rm for}\; \height(\mathfrak{q'})=1.\]
It is enough to show that height of  $\mathfrak{q'} \cap
R/\mathfrak{p}$ is one. Suppose that $\height(\mathfrak{q}) -2$ is
greater than or equal to height of $\mathfrak{p}$. We may assume
that $\mathfrak{q}$ is the maximal ideal of $R$. Then

\[ \height(\mathfrak{q}) -2 \geq \nu(\mathfrak{p})= \ell(\mathfrak{p})=\dim(\Gr/\mathfrak{q}\Gr)=\dim \Gr - \height(\mathfrak{q}\Gr). \]

\noindent Therefore,
 \[\height(\mathfrak{q'}) \geq \height(\mathfrak{q}\Gr) \geq 2.\]

For the converse, suppose that $\mathfrak{p}$ has finite
projective dimension and that $\Gr$ is integrally closed. We make
use of \cite[Theorem 2.4]{JU}, which asserts that the analytic
spread $\ell(\mathfrak{p}_{\mathfrak{q}})$ is less than or equal
to $\height(\mathfrak{q})-2$, for the relevant primes of
$V(\mathfrak{p})$. Since $\mathfrak{p}$ is of linear type, its
analytic spread $\ell(\mathfrak{p}_{\mathfrak{q}})$ and minimum
number of generators $\nu(\mathfrak{p}_{\mathfrak{q}})$ are the
same. \QED

\bigskip

Suppose that the associated graded ring $\Gr$ of $\mathfrak{p}$ is
a domain but not integrally closed. We would like to describe the
set of prime ideals $\mathfrak{q}$ of $R/\mathfrak{p}$ such that
$\Gr_{\mathfrak{q}}$ is not integrally closed. The condition for
the associated graded ring $\Gr$ to be a domain under the
assumptions of Proposition~\ref{appcx1} can be rephrased in terms
of the Fitting ideals of the prime ideal $\mathfrak{p}$. Suppose
that $\mathfrak{p}$ has a presentation
\[R^m \stackrel{\varphi}{\rightarrow} R^n \rightarrow \mathfrak{p}
\rightarrow 0.\] Denote the ideal generated by $t \times t$ minors
of $\varphi$ by $I_t(\varphi)$.  For every prime ideal
$\mathfrak{q}$ which contains $\mathfrak{p}$ properly, the
condition
\[\nu(\mathfrak{p}_{\mathfrak{q}}) \leq
\height(\mathfrak{q}) -1\]is equivalent to
\[\grade I_{t}(\varphi) \geq (n-1)-t+3, \]for $1 \leq t
\leq n- \height(\mathfrak{p})$ (\cite[Corollary 6.7]{HSV3}). For
an integer $s$ such that $\height(\mathfrak{p}) +1 \leq s \leq
\dim R$, letting $t=n-s+1$ gives that
\[\height(I_{n-s+1}(\varphi)) \geq s+1. \]Suppose that
$R/\mathfrak{p}$ is integrally closed. We define the {\em normal
locus} of the associated graded ring $\Gr$ as the set
 \[ \mbox{\rm NL}(\Gr) = \{ \mathfrak{q} \in \spec(R/\mathfrak{p}) \mid  \Gr_{\mathfrak{q}} \ \textrm{is an
integrally closed domain}\}.\] The following proposition shows
that the normal locus of the associated graded ring $\Gr$ is
determined by the Fitting ideals of $\mathfrak{p}$ or equivalently
of $\mathfrak{p}/\mathfrak{p}^2$ (if $\mathfrak{p}$ has sliding
depth).

\begin{proposition}
Let $R$ be a Cohen-Macaulay ring and $\mathfrak{p}$ a prime ideal
with sliding depth. Suppose that $\mathfrak{p}$ satisfies
$G_{\infty}$, that $R/\mathfrak{p}$ is integrally closed and that
the associated graded ring $\Gr$ of $\mathfrak{p}$ is a domain.
Then the set {\rm NL($\Gr$)} is an open subset of
$V(\mathfrak{p})$.
\end{proposition}

\demo We may assume that $\Gr$ is not integrally closed. Suppose
that $\mathfrak{p}$ has a presentation $R^m
\stackrel{\varphi}{\rightarrow} R^n \rightarrow \mathfrak{p}
\rightarrow 0$. Since $\Gr$ is a domain, using
Proposition~\ref{appcx1} and the argument above, we have
\[\height(I_{n-s+1}(\varphi)) \geq
s+1, \] for all $s$ such that $\height(\mathfrak{p}) +1 \leq s
\leq n$.

On the other hand, since the associated graded ring $\Gr$ is not
integrally closed, by Proposition~\ref{appcx2}, there exists a
prime ideal $\mathfrak{q} \supset \mathfrak{p}$ of height $t \geq
\height(\mathfrak{p}) +2$ such that
$\nu(\mathfrak{p}_{\mathfrak{q}} ) > \height(\mathfrak{q})-2$.
Because of the existence of $\mathfrak{q}$,
$\height(I_{n-t+2}(\varphi)) = t$. Now let $\Sigma = \{ t \mid
\height(I_{n-t+2} (\varphi)) = t,\;\;\height(\mathfrak{p}) +2 \leq
t \leq n\}$. For each $t \in \Sigma$, let $I_{n-t+2} (\varphi)=
I_{n-t+2}' \bigcap I_{n-t+2}''$, where $\height(I_{n-t+2}') $ is
$t$ and $\height(I_{n-t+2}'') $ is greater than $t$. Then the set
$\mbox{\rm NL}(\Gr)$ is the complement of $\bigcup_{t \in \Sigma}
V(I_{n-t+2}')$. \QED

\bigskip

Let us put this result under some perspective. Let $S$ be a domain
and $E$ a finitely generated $S$--module. In general, we are not
aware of a way to determine that the symmetric algebra of $E$ is a
domain using only the Fitting ideals of $E$. In contrast, this is
possible for the conormal module with the properties above. The
same observations apply to the obstructions to the normality of
the symmetric algebra of $E$.

\section{Divisor Class Groups}\label{DCG}

Throughout this section, we assume that $\mathfrak{p}$ is a prime
$R$--ideal of finite projective dimension and that the associated
graded ring $\Gr$ of $\mathfrak{p}$ is a domain. Under these
assumptions the associated graded ring of $\mathfrak{p}$ is
isomorphic to the Rees algebra of the conormal
$R/\mathfrak{p}$--module $\mathfrak{p}/\mathfrak{p}^2$. We are
going to compare the divisor class group of $R/\mathfrak{p}$ and
that of $\Gr$ when they are both integrally closed domains. We
recall that if $R$ is an integrally closed domain, then $\Div(R)$
is the group of divisorial ideals with the operation $I \circ J =
((IJ)^{-1})^{-1}$. The {\em divisor class group} $\Cl(R)$ is the
quotient group $\Div(R) / \Prin(R)$, where $\Prin(R)$ is the
subgroup of principal fractional ideals (\cite[Proposition 3.4 and
\S 6]{Fossum}).

\begin{theorem} \label{hong1} Let $R$ be a Cohen--Macaulay ring and $\mathfrak{p}$ a prime ideal of finite projective dimension. If
$R/\mathfrak{p}$ is an integrally closed domain of dimension
greater than or equal to $2$ and the associated graded ring $\Gr$
of $\mathfrak{p}$ is an integrally closed domain, then the mapping
of divisor class groups
\[ \Cl(R/\mathfrak{p})\ni [L] \mapsto [L\Gr] \in \Cl(\Gr),\]
is a group isomorphism. In particular, if $\Gr$ is integrally
closed and $R/\mathfrak{p}$ is a factorial domain, then  $\Gr$ is
a factorial domain.
\end{theorem}

If $R/\mathfrak{p}$ has dimension less than or equal to one, then
$\mathfrak{p}$ is a locally complete intersection. It follows that
$\mathfrak{p}/\mathfrak{p}^2$ is projective and that the
associated graded ring $\Gr$ of $\mathfrak{p}$ is isomorphic to
the symmetric algebra $S(\mathfrak{p}/\mathfrak{p}^2)$, when the
isomorphism between $\Cl(\Gr)$ and $\Cl(R/\mathfrak{p})$ will be
taken care of by Proposition~\ref{flatext}.

Our study of the divisor class group of an associated graded ring
benefits enormously from a result of Johnson \& Ulrich (\cite{JU})
that circumscribe very explicitly the Serre's condition ($R_k$)
for the associated graded ring. In particular, we use
\cite[Theorem 2.4]{JU} in order to prove that the mapping in
Theorem~\ref{hong1} is well defined.

\begin{proposition} \label{preflat}
Let $R$ be a Cohen--Macaulay ring and $\mathfrak{p}$ a prime ideal
of finite projective dimension. Suppose that the associated graded
ring $\Gr$ of $\mathfrak{p}$ is an integrally closed domain. For
$\mathfrak{q}\in \spec(R/\mathfrak{p})$ of height at least $2$,
height of $\mathfrak{q}\Gr$ is greater than or equal to $2$.
\end{proposition}

\demo Suppose that $\mathfrak{q}G$ is contained in a prime ideal
$\mathfrak{P}$ of height $1$. Setting $\mathfrak{m}$ to be the
inverse image of $\mathfrak{P}$ in $R$ and localizing, we may
assume that $(R,\mathfrak{m})$ is a local ring and
$\mathfrak{q}\Gr\subset \mathfrak{m}\Gr\subset \mathfrak{P}$. Now
$\height(\mathfrak{m}\Gr)$ is $1$ and
\[\ell(\mathfrak{p})= \dim \Gr/\mathfrak{m}\Gr = \dim \Gr - \height(\mathfrak{m}\Gr) = \dim R-1,\]
which, by \cite[Theorem 2.4]{JU}, means that

\[ \dim R -1 = \ell(\mathfrak{p}) \leq \max \{\height(\mathfrak{p})\;,\; \dim
R -2 \}=\dim R -2 .\] \QED

\begin{proposition} \label{flatext}
Let $R$ be a Krull domain and $E$ a finitely generated projective
$R$--module. Then the symmetric algebra $S=S_R(E)$ is a Krull
domain and the mapping of divisor class groups
\[ \Cl(R) \longrightarrow \Cl(S), \quad [L] \mapsto [LS],\]
is a group isomorphism.
\end{proposition}

\demo The proof is nearly the same as when the symmetric algebra
$S$ is a ring of polynomials. It consists of two main
observations. Since the morphism $R\rar S$ is flat, there is an
induced homomorphism $\varphi: \Cl(R)\rar \Cl(S)$ of divisor class
groups (\cite[Proposition 6.4]{Fossum}). We claim that $\varphi$
is an isomorphism. Let $I$ and $J$ be two divisorial ideals of
$R$. Note that
\[\Hom_S(I\otimes_R S, J\otimes_RS) = \Hom_R(I,J)\otimes_RS.\]
Let $\varphi=\bigoplus_{i\geq 0} \varphi_i \in \Hom_S(I\otimes_R
S, J\otimes_R S)$ be the isomorphism between $I\otimes_R S$ and $J
\otimes_R S$ with $\phi=\bigoplus_{i \geq 0} \phi_i$ such that
$\varphi \circ \phi = \bigoplus_{i \geq 0} \varphi_i \circ \phi_i$
is the identity map. Since $I$ and $J$ are positively graded, it
means that only $\varphi_0 \circ \phi_0$ is the identity map and
the others must vanish. Therefore any isomorphism between
$I\otimes_R S$ and $J \otimes_R S$ must be realized by an
isomorphism between $I$ and $J$, which shows that $\varphi$ is
injective.

In order to show that $\varphi$ is surjective, it suffices to
prove that for every divisorial prime ideal $\mathfrak{P}$ of $S$,
its divisor class $[\mathfrak{P}]$ lies in the image of $\varphi$.
Suppose that $\mathfrak{P}\cap R = \mathfrak{q}\neq 0$. Since
$S_R(E)/\mathfrak{q}S_R(E) \simeq S_{R/\mathfrak{q}}
(E/\mathfrak{q}E)$ which is a domain, $\mathfrak{q}S$ is a prime
ideal so that $\mathfrak{P}=\mathfrak{q}S$. Suppose that
$\mathfrak{q}=0$. Let $0\neq f\in \mathfrak{P}$ be an element that
generates the extension of $\mathfrak{P}$ to $KS$, where $K$ is
the field of
fractions of $R$. 
We recall that $S$ is a projective $R$--module and that $f\in Q
\subset S$, where $Q$ is a direct summand of $S$. Let $I$ be the
$R$--ideal
\[ I = (\; \alpha(f) \mid \alpha \in \Hom_R(Q,R)).\]
Then $\mathfrak{P}= I^{-1}fS$, 
which implies that $\varphi([I^{-1}])= [I^{-1}S]=[\mathfrak{P}]$.
\QED

\bigskip

We are ready to complete the proof of Theorem~\ref{hong1}.

\bigskip

\noindent {\bf  Proof of Theorem~\ref{hong1}:} For any prime ideal
of $R/\mathfrak{p}$ of the form $\mathfrak{q}/\mathfrak{p}$, of
height at least $2$, there is a regular sequence $\alpha, \beta $
in $\mathfrak{q}/\mathfrak{p}$, and therefore
$(\mathfrak{q}/\mathfrak{p})[x]$ contains the prime element
$f=\alpha+\beta x$. We claim that $f$ is also a prime element in
$\Gr[x]$. This follows from a series of observations. First, any
minimal prime of $(\alpha,\beta)$ in $R/\mathfrak{p}$ has height
$2$ and therefore by Proposition~\ref{preflat}, the extension
$(\alpha,\beta)\Gr$ has height $2$. Since $\Gr$ is integrally
closed, the elements $\alpha$ and $\beta$ generate an $\Gr$--ideal
of grade two and $f=\alpha + \beta x$ is a prime element of
$\Gr[x]$ (\cite[Lemma 14.1]{Fossum}).

Another elementary property of the calculation of divisor class
groups is that it is unaffected when the rings are localized with
respect to multiplicative sets formed by powers of prime elements.
Let $W$ be a multiplicative set generated by \[\{ \alpha + \beta x
\mid \alpha , \beta \;\;\mbox{is}
\;\;R/\mathfrak{p}\mbox{--regular sequence} \}.\] The elements
$\alpha + \beta x$ in $W$ are prime elements in both
$(R/\mathfrak{p})[x]$ and $\Gr[x]$. Consider the following
commutative diagram of divisor class groups.

\bigskip

\[
\begin{CD}
\Cl(R/\mathfrak{p})       @ >>> \Cl(\Gr) \\
@V{\cong}VV       @VV{\cong}V\\
\Cl((R/\mathfrak{p})[x])    @ >>> \Cl(\Gr[x])\\
@V{\cong}VV        @VV{\cong}V\\
\Cl((R/\mathfrak{p})[x]_W ) @>>> \Cl(\Gr[x]_W )
\end{CD}
\]

\bigskip

\noindent Denoting $(R/\mathfrak{p})[x]_W$ by $A/J$, the ideal $J$
must be locally a complete intersection by \cite[Theorem 2.4]{JU} and
\cite{CowsikNori}.
 This means that $J/J^2$ is a finitely
generated projective $A/J$--module and that $\Gr[x]_W$ is
isomorphic to $S_{A/J}(J/J^2)$ (\cite[Theorem 6.1]{HSV3}). By
Proposition~\ref{flatext}, we complete the proof. \QED

\bigskip

Now we explore the homomorphism between divisor class groups in
the case where $R/\mathfrak{p}$ is integrally closed but the
associated graded ring $\Gr$ of $\mathfrak{p}$ is a domain that is
not integrally closed. We assume that the integral closure
$\overline{\Gr}$ of $\Gr$ is finitely generated as an
$\Gr$--module (a condition that holds true when $R$ is an affine
ring over a field, and in a greater generality).

\begin{theorem} \label{hong2} Let $R$ be a Cohen--Macaulay ring and
$\mathfrak{p}$ a prime ideal of finite projective dimension.
Suppose that $S=R/\mathfrak{p}$ is an integrally closed domain and
that the associated graded ring $\Gr$ of $\mathfrak{p}$ is a
domain with finite integral closure $\overline{\Gr}$. Then there
exists an exact sequence of divisor class groups
\[ 0 \rar H \lar \Cl(\overline{\Gr}) \lar \Cl(R/\mathfrak{p})
\rar 0,\] where $H$ is a finitely generated subgroup of
$\Cl(\overline{\Gr})$ that vanishes if $\Gr$ is integrally closed.
\end{theorem}

\demo We make two observations. First, for each prime ideal
$\mathfrak{q}=Q/\mathfrak{p}$ of $S=R/\mathfrak{p}$ of height one,
localizing at $Q$ we get $QR_Q=(\mathfrak{p},x)R_Q$, since
$S_{\mathfrak{q}}$ is a discrete valuation ring. From the exact
sequence
\[ 0 \rar R_Q/\mathfrak{p}R_Q \stackrel{x}{\lar} R_Q/\mathfrak{p}R_Q \lar R_Q/QR_Q \rar 0,\]
it follows that the residue field of $R_Q$ has finite projective
dimension since $R/\mathfrak{p}$ has finite projective dimension
by hypothesis. Therefore the local ring $(R_Q,QR_Q)$ is a regular
local ring by Serre's theorem \cite[Theorem 2.2.7]{BH}. As a
consequence $\mathfrak{p}R_Q$ is a complete intersection. We have
thus shown that $\mathfrak{p}$ is a complete intersection in
$\height(\mathfrak{p}) + 1$.

Next, we may assume that the associated graded ring $\Gr$ is not
integrally closed. Let $C= \ann(\overline{\Gr}/\Gr)$ be the
conductor ideal. Since $\Gr$ and $\overline{\Gr}$ are graded
rings, $C$ is a homogeneous ideal of $\Gr$. Set
\[ L = C\cap S,\]
the component of $C$ in degree zero. We claim that $L$ is not
contained in any prime ideal of $S$ of height one. Suppose that
$L$ is contained in a prime ideal $\mathfrak{q}=Q/\mathfrak{p}$ of
height one. Then we have $L_{\mathfrak{q}}\subset
\mathfrak{q}S_{\mathfrak{q}}$. But $\mathfrak{p}R_Q$ is a complete
intersection and therefore
$\Gr_{\mathfrak{q}}=\gr_{\mathfrak{p}R_Q} R_Q$ is a ring of
polynomials over the discrete valuation ring $S_{\mathfrak{q}}$.
This means that $\Gr_{\mathfrak{q}}$ is
$\overline{\Gr}_{\mathfrak{q}}$ and $C_{\mathfrak{q}}$ is
$S_{\mathfrak{q}}$.

By this observation, the $S$--ideal $L$ has height greater than or
equal to $2$. Since $S$ is integrally closed, there is a
$S$-regular sequence $a,b$ contained in $L$. Consider the addition
of an indeterminate $x$ to $S$, $\Gr$ and $\overline{\Gr}$. We
consider the prime element $f=a+bx$ of $S[x]$. Since $f\in
\ann(\overline{\Gr}[x]/\Gr[x])$, we get $\overline{\Gr}[x]_f=
\Gr[x]_f$. But $\Gr[x]_f$ is just the associated graded ring of
the ideal $\mathfrak{p}R[x]_f$. This ideal inherits all the
properties of $\mathfrak{p}$, in particular it is a prime ideal of
finite projective dimension. Applying Theorem~\ref{hong1}, we
obtain these isomorphisms of divisor class groups
\[ \Cl( \overline{\Gr}[x]_f) = \Cl(\Gr[x]_f ) \cong \Cl(\gr_{\mathfrak{p}R[x]_f} R[x]_f )
\cong \Cl(R[x]_f / \mathfrak{p}R[x]_f ) \cong \Cl(R/\mathfrak{p}).
\] In particular, the last isomorphism follows from the fact that $f$ is a prime element of $S[x]$.
On the other hand, using the exact sequence associated with the
localization formula for divisor class groups \cite[Corollary
7.2]{Fossum}, we obtain
\[ 0 \rar H \lar \Cl(\overline{\Gr}[x]) \lar \Cl(
\overline{\Gr}[x]_f) \rar 0,\] where $H$ is generated by the
classes of all primes divisors in $\overline{\Gr}[x]$ that contain
$f$ so that $H$ is finitely generated \cite[Proposition
1.9]{Fossum}. Finally, replacing $\Cl(\overline{\Gr}[x]_f)$ with
$\Cl(R/\mathfrak{p})$ and $\Cl(\overline{\Gr}[x]) $ with
$\Cl(\overline{\Gr}) $ gives the desired sequence. \QED

\bigskip

While we do not know the structure of $H$ in detail, examples and
general arguments suggest that the following holds.

\begin{conjecture}\label{conj1}
{\rm $H$ is a free group.}
\end{conjecture}

\begin{example}{\rm Let $R$ be the polynomial ring $k[u,v,t,w]$, $\mathfrak{p}$
the $R$--ideal defining $S=R/\mathfrak{p}=k[x^3 , x^2 y, xy^2,
y^3]$, and $\mathfrak{m}$ the $S$-ideal $(x^3 , x^2 y, xy^2,
y^3)$. Using Macaulay2, we compute the associated graded ring
$\Gr$ and the integral closure $\overline{\Gr}$ of $\Gr$.
\[
\begin{array}{lll}
\Gr &= &R[x_1, x_2, x_3]/(\mathfrak{p}, -tx_1+vx_2-ux_3,
wx_1-tx_2+vx_3)\\
\overline{\Gr}&= &\Gr[Y]/(Yw-tx_3, Yt-vx_3, Yv-ux_3,
Yu+vx_1-ux_2,Y^2-Yx_2+x_1x_3).
\end{array}
\]
Then $H$ is generated by one element
$[\mathfrak{m}\overline{\Gr}]$ because
\[
\overline{\Gr}/\mathfrak{m}\overline{\Gr} =
k[x_1,x_2,x_3,Y]/(Y^2-Yx_2+x_1x_3),
\]where the monic quadratic polynomial is irreducible.
We claim that $[\mathfrak{m}\overline{\Gr}]$ is a torsionfree
element. Suppose that $[\mathfrak{m}^n \overline{\Gr}]$ is equal
to $[f\overline{\Gr}]$ for some $n$, which means that
$(\mathfrak{m}^n \overline{\Gr})^{-1-1}=f\overline{\Gr}$. Since
$\mathfrak{m}^n \overline{\Gr}$ is a homogeneous ideal of height
1, $f$ is not a unit and moreover has degree $0$. Now
$\mathfrak{m}^n$ is contained in $f \overline{\Gr}$ so that it is
contained in $fS$, which is impossible because
$\height(\mathfrak{m}^n)$ is $2$ but $\height(fS)$ is $1$. This
proves that $H$ is isomorphic to $\mathbb{Z}$. Moreover,
\cite[Proposition 11.4]{Fossum} shows that $\Cl(S)$ is
$\mathbb{Z}$. Therefore, $\Cl(\overline{\Gr})$ is $\mathbb{Z}
\bigoplus \mathbb{Z}$.}\end{example}

Although we have concentrated on the Rees algebra
$\Rees(\mathfrak{p}/\mathfrak{p}^2)$, there are related algebras
to which the techniques employed here may be applied. For example,
suppose that $R_{\mathfrak{p}}$ is a regular local ring and that
the $\mathfrak{p}$--symbolic filtration $\{ \mathfrak{p}^{(n)}, n
\geq 0 \}$ equals the filtration $\{\overline{\mathfrak{p}^n}, n
\geq 0 \}$. Then the reduced ring $\Gr_{\mbox{\rm red}}$ of the
associated graded ring $\Gr$ is also a domain (\cite[Theorem
2.1]{Huneke1}). The significant difference between the algebras
$\Gr$ and $\Gr_{\mbox{\rm \ red}}$ is that we have criteria for
the ($R_1$)--condition for the associated graded ring $\Gr$ only,
which is essential for the study of normality. Nevertheless we are
able to obtain a similar result to Theorem~\ref{hong2} in case of
the integral closure $\overline{\Gr_{\mbox{\rm red}}}$ of the
reduced ring $\Gr_{\mbox{\rm red}}$.

\begin{proposition}\label{hong7} Let $R$ be a Cohen--Macaulay ring,
$\mathfrak{p}$ a prime ideal of finite projective dimension, and
$\Gr$ the associated graded ring of $\mathfrak{p}$. Suppose that
$S=R/\mathfrak{p}$ is an integrally closed domain and that the
reduced ring $\Gr_{\mbox{\rm  red}}$ of $\Gr$ is a domain
with finite integral closure $\overline{\Gr_{\mbox{\rm
red}}}$. Then there exists an exact sequence of divisor class
groups
\[ 0 \rar H \lar \Cl(\overline{\Gr_{\mbox{\rm red}}}) \lar \Cl(R/\mathfrak{p})
\rar 0,\] where $H$ is a finitely generated subgroup of
$\Cl(\overline{\Gr_{\mbox{\rm red}}})$.
\end{proposition}

\demo If $\Gr_{\mbox{\rm  red}}$ is integrally closed, then $\Gr$
is a domain (\cite[Proposition 2.2]{Huneke1}) so that
$\Gr_{\mbox{\rm   red}}$ is just $\Gr$. Now we may assume that
$\Gr_{\mbox{\rm   red}}$ is not integrally closed. Let $C=
\ann(\overline{\Gr_{\mbox{\rm   red}}}/\Gr_{\mbox{\rm red}})$ be
the conductor ideal. Then the the component of $C$ in degree zero
$L = C\cap S$ has grade at least two. We consider the prime
element $f=a+bx$ of $S[x]$, where $a,b$ is a $S$--regular sequence
contained in $L$ and $x$ is an indeterminate.

Since $f\in \ann(\overline{\Gr_{\mbox{\rm
red}}}[x]/\Gr_{\mbox{\rm   red}}[x])$, we get
$\overline{\Gr_{\mbox{\rm   red}}}[x]_f= \Gr_{\mbox{\rm
red}}[x]_f =\Gr[x]_f$, where the last equality is again from the
result \cite[Proposition 2.2]{Huneke1}. Now we apply
Theorem~\ref{hong1} and obtain the following isomorphisms of
divisor class groups.

\[ \Cl( \overline{\Gr_{\mbox{\rm   red}}}[x]_f) = \Cl(\Gr[x]_f ) \cong \Cl(\gr_{\mathfrak{p}R[x]_f} R[x]_f )
\cong \Cl(R[x]_f / \mathfrak{p}R[x]_f ) \cong \Cl(R/\mathfrak{p}).
\] By \cite[Corollary 7.2]{Fossum}, there is an exact sequence
\[ 0 \rar H \lar \Cl(\overline{ \Gr_{\mbox{\rm   red}}}[x]) \lar \Cl(
\overline{\Gr_{\mbox{\rm   red}}}[x]_f) \rar 0.\] Finally, by
replacing $\Cl(\overline{\Gr_{\mbox{\rm   red}}}[x]_f)$ with
$\Cl(R/\mathfrak{p})$ and $\Cl(\overline{\Gr_{\mbox{\rm
red}}}[x]) $ with $\Cl(\overline{\Gr_{\mbox{\rm   red}}})$, we
show that there is a surjective group homomorphism from
$\Cl(\overline{\Gr_{\mbox{\rm   red}}})$ to
$\Cl(R/\mathfrak{p})$. \QED

\section{Integrally Closed Conormal Modules}\label{ICCM}

If $R$ is an integrally closed domain and $E$ is a torsionfree
$R$--module, the integral closure of the Rees algebra
$\Rees(E)=\bigoplus_{n \geq 0} E_n$ of the module $E$ is the
algebra $\overline{\Rees(E)}=\bigoplus_{n \geq 0} \overline{E_n}$,
where $\overline{E_n}$ is the integral closure of the module $E_n$
for all $n$ (See Proposition~\ref{intclosrees} below). Partly for
this reason, it is worthwhile to study the integral closure of a
module and apply its techniques to the components of the
associated graded ring. Throughout this section we assume that
height of $\mathfrak{p}$ is at least two. If $\mathfrak{p}$ is a
prime $R$--ideal which is generically a complete intersection and
whose associated graded ring $\Gr$ is a domain, then $\Gr$ is
isomorphic to the Rees algebra
$\Rees(\mathfrak{p}/\mathfrak{p}^2)$ of the conormal module
$\mathfrak{p}/\mathfrak{p}^2$. Under such assumptions, we examine
the relationship between the completeness of the components of
$\Gr$ and the normality of $\Gr$.

\begin{proposition}\label{intclosrees} Let $R$ be a normal domain and $E$
a torsionfree finitely generated  $R$--module. The Rees algebra
$\Rees(E) = \bigoplus_{n \geq 0} E_n$ of the module $E$ is
integrally closed if and only if $E_n$ is integrally closed for
all $n$.
\end{proposition}

\demo If each component $E_n$ is integrally  closed, we have
\[ \Rees(E)=\sum_{n\geq 0} E_n =\sum_{n\geq 0}\bigcap_{V}VE_n =  \bigcap_{V}V\Rees(E) = \bigcap_{V}\Rees(VE),\]
where $V$ runs over all the valuation overrings of $R$. Since $VE$
is a free $V$--module, the Rees algebra $\Rees(VE)$ is a ring of
polynomials over $V$. This gives a representation of $\Rees(E)$ as
an intersection of polynomial rings and it is thus normal. The
converse is similar. \QED

\bigskip

Suppose that the associated graded ring $\Gr=\bigoplus_{t \geq 0}
G_t$ satisfies the ($S_2$)--condition. Denote the integral closure
of $\Gr$ by $\overline{\Gr}$ and the ideal $\bigoplus_{t \geq 1}
G_t$ by $\Gr_{+}$. Let $I$ be the ideal generated by
$\Gr$--regular sequence of length two in $\Gr_{+}$. Then
$\depth(I,\overline{\Gr}/\Gr)$ is at least one and hence we can
choose a $\overline{\Gr}/\Gr$--regular element $x$ from the degree
one component. Suppose that $\Gr_{t_0}$ is not integrally closed
for some $t_0$. The map defined as the multiplication by $x$ gives
the embedding
\begin{eqnarray*}\label{gt0} \left(\overline{\Gr}/\Gr \right)_{t_0} \stackrel{\cdot x}{\longrightarrow}
\left(\overline{\Gr}/\Gr \right)_{t_0 +1}, \end{eqnarray*} and
$\Gr_t$ is not integrally closed for every $t \geq t_0$. We shall
be concerned with the following broad conjecture.

\begin{conjecture} \label{conj2} {\rm Let $R$ be a Cohen--Macaulay ring and $\mathfrak{p}$ a
prime ideal of finite projective dimension. Suppose that the
associated graded ring $\Gr$ is a Cohen-Macaulay domain and that
$R/\mathfrak{p}$ is integrally closed. If $\Gr_t=
\mathfrak{p}^t/\mathfrak{p}^{t+1}$ is integrally closed for
$t=\dim R/\mathfrak{p} -1$, then $\Gr$ is normal.
}\end{conjecture}

If $R/\mathfrak{p}$ has dimension less than or equal to one, then
the associated graded ring $\Gr$ is isomorphic to the symmetric
algebra of the projective $R/\mathfrak{p}$--module
$\mathfrak{p}/\mathfrak{p}^2$. Thus far we have settled the
following case.

\begin{theorem} \label{hong3} Let $R$ be a Gorenstein local ring and $\mathfrak{p}$ a prime ideal generated by a strongly Cohen-Macaulay
$d$--sequence. Suppose that $\mathfrak{p}$ has finite projective
dimension, that $S=R/\mathfrak{p}$ is an integrally closed domain
of dimension $2$, and that the associated graded ring $\Gr$ of
$\mathfrak{p}$ is a domain. Then the conormal module
$E=\mathfrak{p}/\mathfrak{p}^2$ is integrally closed if and only
if $\Gr$ is normal.
\end{theorem}

Before we prove Theorem~\ref{hong3}, we briefly recall how to
attach divisors to certain modules (see \cite{Fossum}). Let $R$ be
a Noetherian normal domain and $E$ a finitely generated
torsionfree $R$--module of rank $r$. The {\em determinant} of $E$
is rank one module $ \det(E) = \wedge^r E$, while the {\em
determinantal divisor}
 of $E$ is the divisor class $\div(E)$
of the bidual $(\det(E))^{-1-1}$. The following are two elementary
properties of this construction (we recall that the codimension
$\codim A$ of a module $A$ is the codimension of its annihilator).
The first of these has an immediate proof.

\begin{proposition}\label{diva}
Let $R$ be a Noetherian normal domain  and let
\[ 0 \rar A \lar B \lar C \lar D \rar 0,\] be an exact sequence of
finitely generated $R$--modules.  If $\codim A\geq 1$ and $\codim
D\geq 2$, then $\div(B)=\div(C)$. In particular, for any finitely
generated torsionfree $R$--module $E$, $\div(E)= \div(E^{-1-1})$.
\end{proposition}

\begin{proposition}\label{div3} Let $R$ be a Noetherian normal domain. Suppose that the complex of finitely generated
$R$--modules
\[ 0 \rar A \stackrel{\varphi}{\lar} B \stackrel{\psi}{\lar} C \rar 0\]
 is an exact sequence
of free modules in each localization $R_{\mathfrak{p}}$ at height
one prime ideals. Then $ \div(B)= \div(A) + \div(C) $. In
particular, if $B$ is a free $R$--module, then $\div(A)=-\div(C)$.
\end{proposition}

\demo We break up the complex into  simpler exact  complexes:
\[ 0 \rar \ker(\varphi) \lar A \lar \im(\varphi)  \rar 0, \]

\[ 0 \rar \im(\varphi) \lar \ker(\psi) \lar \ker(\psi)/\im(\varphi) \rar 0, \]

\[0 \rar \im(\psi) \lar C \lar C/\im(\psi) \rar 0,\]
By hypothesis, $\codim \ker(\varphi)\geq 2$, $\codim
\ker(\psi)/\im(\varphi) \geq 2$, $\codim\; C/\im(\psi) \geq 2$. By
Proposition~\ref{diva} we have the equality of determinantal
divisors.

\[ \div(A)=\div(\im(\varphi))= \div(\ker(\psi)),\;\; \div(C)=
\div(\im(\psi)).\] What this all means is that we may assume the
given complex is exact.

\smallskip

Suppose $\rank(A)=r$ and $\rank(C)=s$ and set $n=r+s$. Consider
the pair $\wedge^rA$, $\wedge^sC$. For $v_1, \ldots, v_r\in A$,
$u_1, \ldots, u_s \in C$, pick $w_i$ in $B$ such that
$\psi(w_i)=u_i$ and consider
\[ v_1 \wedge \cdots \wedge v_r \wedge w_1 \wedge \cdots \wedge
w_s\in \wedge^{n}B. \]
Different choices for $w_i$ would produce elements in $\wedge^{n}B$
that differ from the above by terms that would contain at least $r+1$
factors of the form
\[ v_1 \wedge \cdots \wedge v_r \wedge v_{r+1} \wedge \cdots,\]
with $v_i \in A$. Such products are torsion elements in $\wedge^nB$.
This implies that modulo torsion we have a well defined pairing
\[ [\wedge^r A/\textrm{\rm torsion}] \otimes_R [\wedge^s C/\textrm{\rm
torsion}] \lar [\wedge^n B/\textrm{\rm torsion}].
\] When localized at primes $\mathfrak{p}$ of codimension at most $1$
the complex becomes an exact complex of projective
$R_{\mathfrak{p}}$--modules and the pairing is an isomorphism. Upon
taking biduals and the $\circ $ divisorial composition, we obtain the asserted
isomorphism. \QED

\bigskip

\noindent{\bf Proof of Theorem~\ref{hong3}:} Let $n$ be the
minimal number of generators of $\mathfrak{p}$. Since the
associated graded ring $\Gr$ is a domain, $n$ is either $\dim R-2$
or $\dim R -1$ (Proposition~\ref{appcx1}). Suppose that the
associated graded ring $\Gr$ is not integrally closed and that the
conormal module $E=\mathfrak{p}/\mathfrak{p}^2$ is integrally
closed. If $n$ is equal to $\dim R-2$, then $\mathfrak{p}$ is
generated by a regular sequence and so we may assume that $n$ is
equal to $\dim R-1$. Since $R_{\mathfrak{p}}$ is a regular local
ring, $\rank E$ is equal to $\height(\mathfrak{p})$, that is
$n-1$. Let
\[E=\overline{E} \hookrightarrow E_0 \hookrightarrow S^{n-1}\]
be the embedding of $E$ into its bidual $E_0=\Hom_S ( \Hom_S (E,S)
, S) $.

At this point we need the notion of $\mathfrak{m}$--full modules.
A torsionfree $R$--module $E$ is called an {\em
$\mathfrak{m}$--full module} if there is an element $x \in
\mathfrak{m}$ such that $\mathfrak{m}E:_{R^r} x =E$. Integrally
closed modules are $\mathfrak{m}$--full modules (\cite[Proposition
2.6]{BV}) and since $\lambda(E_0 /E)$ is finite, the minimal
number $\nu (E_0 )$ of generators of $E_0$ is either $n$ or $n-1$
(\cite[Corollary 2.7]{BV}).

Suppose $\nu (E_0 )$ is equal to $n-1$. Consider the following
diagram.

\[
\begin{CD}
0 \rightarrow H_1=\kernel (\phi)     @>>> S^n      @>{\phi}>>    E @>>> 0 \\
              \;\;\;\;\;\;@V{=}VV               @V{=}VV                @VV{i}V  \\
0 \rightarrow \kernel (\varphi)   @>>> S^n      @>>{\varphi}> E_0
=S^{n-1} @>>> E_0/E\rightarrow 0
\end{CD}
\]

\noindent where the short exact sequence in the top row is from
the approximation complex of $E$ and $i$ is an inclusion such that
$\varphi = i \circ \phi$. In the exact sequence in the bottom row,
by using Proposition~\ref{div3}, we get $\div(\kernel (\varphi)) =
\div(E_0/E)$, and hence

\[S \simeq \kernel(\varphi) =\kernel(\phi) = H_1.\]
This gives a contradiction because $\mathfrak{p}$ has finite
projective dimension (\cite[Theorem 1.4.9]{G},\cite[Theorem
5.2.1]{V1}).

Suppose $\nu (E_0 )$ is equal to $n$. Consider the following
diagram.

\[
\begin{CD}
0 \longrightarrow  K  @>>>  S^n  @>d_0>>  E_0  @>>> 0 \\
 \;\;\;\;\;\;\;\;\;@AAA @A{A}AA                 @AA{i}A  \\
 0  \longrightarrow H_1  @>>>  S^n @>>{d}> E  @>>> 0
\end{CD}
\]We claim that $K$ is isomorphic to the first Koszul homology module $H_1$ associated to $\mathfrak{p}$. Since $E_{\mathfrak{q}}$
is reflexive for every height $1$ prime $\mathfrak{q}$, the
determinantal divisor $\div(E)$ of $E$ equals $\div(E_0)$. By
using Proposition~\ref{div3}, we obtain that $\div(K)$ equals
$\div(H_1)$. Since $K$ and $H_1$ are divisorial, the claim is
proved. Let $a,b$ be a $S$-regular sequence. By tensoring the
exact sequence in the top row with $S/(a,b)S$, we have the short
exact sequence

\[ 0 \rightarrow \overline{K} \rightarrow \overline{S^n} \rightarrow
\overline{E_0} \rightarrow 0 \;,\]

\noindent which splits because $\overline{H_1}$ is injective.
Therefore, $\overline{E_0}$ is isomorphic to $\overline{S^{n-1}}$,
which is a contradiction. \QED

\bigskip

Although this is far less than what we would wish, it will suffice
to develop methods to set up the computation of the integral
closure of some associated graded domains. Under the conditions of
the theorem, we actually prove the following corollary.

\begin{corollary}\label{ref1} If $S= R/\mathfrak{p}$ is an integrally closed domain, of arbitrary
dimension, then the conormal module
$E=\mathfrak{p}/\mathfrak{p}^2$ is integrally closed if and only
if it is reflexive.
\end{corollary}

\demo Consider the following approximation complex
\[ 0 \rightarrow H_1 \rightarrow S^n \rightarrow E
\rightarrow 0.\] Using the Depth Lemma(\cite[Lemma 3.1.4]{V2}),
\[\depth (E) \geq \min \{ \depth (H_1) -1 , \depth (S^n) \}= \dim
S -1.
\] Now it is enough to consider prime $S$--ideals whose height is less than or equal to two.
For every prime ideal $\mathfrak{q}$ of height one,
$E_{\mathfrak{q}}$ is free because it is a finitely generated
torsionfree module over a discrete valuation ring. For prime
ideals $\mathfrak{q}$ of height two, Theorem~\ref{hong3} shows
that $E_{\mathfrak{q}}$ is free. Therefore $E$ satisfies the
$(S_2)$--condition. \QED

\bigskip

Now we give an independent proof of Theorem~\ref{hong3} when the
prime ideal $\mathfrak{p}$ has height $2$. It indicates where we
should look the integral closure of the conormal module.
 It is a consequence of
very useful criterion of completeness.

\begin{proposition}\label{div2}
Let $R$ be a Cohen--Macaulay local ring and $\mathfrak{p}$ a prime
$R$--ideal of height $2$ with a finite free resolution
\[
0 \rightarrow R^2 \stackrel{\phi}{\rightarrow} R^3 \rightarrow
\mathfrak{p} \rightarrow 0.
\]
Suppose that $S=R/\mathfrak{p}$ is integrally closed. Then the
determinant $\det(E)$ of the conormal module
$E=\mathfrak{p}/\mathfrak{p}^2$ is divisorial.
\end{proposition}

\demo The conormal module $E$ has rank $2$ with an embedding $E
\hookrightarrow S^2$. Composing the embedding with the free
resolution of $E$, we have the following exact sequence

\[
S^2 \stackrel{\overline{\phi}}{\longrightarrow} S^3
\stackrel{\varphi}{\longrightarrow} S^2,\]
\[\phi= \left[\begin{array}{ll} a & d \\
b &e\\c&f \end{array} \right]\;,\;\mbox{and}\;\; \varphi= \left[\begin{array}{lll} a_1 & b_1 & c_1 \\
a_2 & b_2 & c_2 \end{array} \right],
\]
where the image of $\varphi$ is $E$. The ideal $J=(a,b,c)$ of $R$
contains $\mathfrak{p}$ properly so that $J$ is a complete
intersection. Since the $S$--ideal $J'=J/\mathfrak{p}$ has no
embedded primes, $J'$ is divisorial. Using $0=\varphi \circ
\overline{\phi}$, we have

\[
\frac{b_1 c_2 - b_2 c_1}{\overline{a}} = \frac{a_2 c_1 - a_1
c_2}{\overline{b}} = \frac{a_1 b_2 - a_2 b_1}{\overline{c}}.
\]

\noindent Then $\det(E)=I_2 (\varphi)=\delta J'$, which proves
that $\det(E)$ is divisorial.\QED

\begin{proposition}\label{div1}
Let $S$ be an integrally closed domain and $E$ a finitely
generated integrally closed torsionfree $S$-module whose order
determinant $\det (E)$ is divisorial. Then $E$ is reflexive.
\end{proposition}

\demo Let $E_0=\Hom_S(\Hom_S(E,S),S)$ and $x \in E_0 \setminus E$.
Since $(E_0)_{\mathfrak{p}}=E_{\mathfrak{p}}$ for every prime
ideal $\mathfrak{p}$ of height one, the ideal $L=E:x$ has height
greater than one. Choose elements $e_1, \ldots, e_{r-1}$ from $E$.
Then

\[L(x \wedge e_1 \wedge \ldots \wedge e_{r-1})=Lx \wedge e_1 \wedge \ldots \wedge e_{r-1} \subset \det (E)\]
Since $\height(L) \geq 2$ and $\det(E)$ is divisorial, $x \wedge
e_1 \wedge \ldots \wedge e_{r-1} \in \det(E)$, i.e.,
$\det(x,E)=\det(E)$. Recall that if $\det(F)=\det(E)$ for $F
\supset E$, then $F \subseteq \overline{E}$. Now $(x,E)$ is
integral over $E$, which means that $x \in \overline{E}=E$. \QED

\begin{example}\label{ex1}{\rm
Let $R$ be the polynomial ring $k[u,v,t,w]$ and $\mathfrak{p}$ the
$R$-ideal generated by $(uw-tv, ut-v^2, vw-t^2 )$. Note that
$S=R/\mathfrak{p}$ is isomorphic to $k[x^3, x^2 y, xy^2, y^3 ]$.
Let $E$ be the conormal module $\mathfrak{p}/\mathfrak{p}^2$,
$E_0$ its bidual $\Hom_S(\Hom_S(E,S),S)$. It is easy to check that
the ideal $\mathfrak{p}$ satisfies the assumptions of
Theorem~\ref{hong3}. In particular the $R$-ideal $\mathfrak{p}$
and the $S$-module $E$ have the following presentations

\[
\begin{array}{ll}
0 \longrightarrow R^2 \stackrel{\phi}{\longrightarrow} R^3
\longrightarrow \mathfrak{p} \longrightarrow 0,&\phi = \left[
\begin{array}{rr} -t & w \\ v & -t \\ -u & v
\end{array} \right]\;, \\
&  \\ S^2 \stackrel{\overline{\phi}}{\longrightarrow} S^3
\stackrel{\varphi}{\longrightarrow} S^2,&\varphi = \left[
\begin{array}{rrr} w &0&-t \\ 0&u & -v
\end{array} \right]\;,
\end{array}
 \]

\noindent where the image of $\varphi$ is $E$. Using Macaulay2, we
compute the associated graded ring $\Gr$ and its integral closure
$\overline{\Gr}\;$.

\bigskip

\[
\begin{array}{lll}
\Gr &= &R[x_1, x_2, x_3]/(\mathfrak{p}, -tx_1+vx_2-ux_3,
wx_1-tx_2+vx_3),\\
\overline{\Gr}&= &\Gr[Y]/(Yw-tx_3, Yt-vx_3, Yv-ux_3,
Yu+vx_1-ux_2,Y^2-Yx_2+x_1x_3).
\end{array}
\]

\bigskip

\noindent Since the module $E$ is not reflexive, $E$ is not
integrally closed (Corollary~\ref{ref1}). Having
$\det(E)=\det(E_0)=v(t,w,u)$ gives that the integral closure
$\overline{E}$ is the bidual module $E_0$. Furthermore we claim
that $\overline{\Gr}$ is Cohen--Macaulay. Let $L$ be
$\Gr[Y]$--ideal generated by $(Yw-tx_3, Yt-vx_3, Yv-ux_3,
Yu+vx_1-ux_2)$. Note that $Y^2-Yx_2+x_1x_3$ is a monic irreducible
polynomial. From the exact sequence
\[ 0 \rightarrow L \rightarrow \Gr[Y]/(Y^2-Yx_2+x_1x_3)
\rightarrow \overline{\Gr} \rightarrow 0,\]we obtain
\[\Hom_{\Gr}(\overline{\Gr},\Gr)=(\Gr :_{\Gr} \overline{\Gr}) \simeq (\beta \mid \alpha + \beta Y \in L)=\mathfrak{m}\Gr. \]
Since $\Gr$ is Gorenstein, by \cite[Proposition 1.1.15]{V2},
$\overline{\Gr}$ is reflexive, and hence
\[ \overline{\Gr} = \Hom_{\Gr}(\Hom_{\Gr}(\overline{\Gr}, \Gr),\Gr) \simeq \Hom_{\Gr}(\mathfrak{m}\Gr, \Gr).\]
Since $\mathfrak{m}\Gr$ is Cohen--Macaulay,
the integral closure $\overline{\Gr}$ is Cohen--Macaulay. 
}
\end{example}

\bigskip

We are going to study the module theoretic properties of some
components of $\Gr$. It is helpful to find the integral closure in
some cases.

\begin{proposition}\label{ddual1}
Let $R$ be a Gorenstein local ring, $\mathfrak{p}$ a perfect prime
ideal generated by a strongly Cohen--Macaulay $d$--sequence, and
$\Gr=\bigoplus_{i \geq 0} \Gr_i$ the associated graded ring of
$\mathfrak{p}$. Suppose that $S=R/\mathfrak{p}$ is an integrally
closed domain and that for every proper prime ideal $\mathfrak{q}
\supset \mathfrak{p}$,
\[\nu(\mathfrak{p}_{\mathfrak{q}})\leq \height(\mathfrak{q})-1,\]
where $\nu(\mathfrak{p}_{\mathfrak{q}})$ is the minimal number of
generators of $\mathfrak{p}_{\mathfrak{q}}$. Denote the minimal
number of generators of $\mathfrak{p}$ by $n$ and height of
$\mathfrak{p}$ by $g$. Then $\Gr_{n-g}$ is not reflexive.
\end{proposition}

\demo Suppose that $\Gr_{n-g}$ is reflexive. By
Proposition~\ref{appcx1}, the associated graded ring $\Gr$ is a
domain. By applying the Depth lemma(\cite[Lemma 3.1.4]{V2}) to the
approximation complex

\[
0 \rightarrow H_{n-g} \stackrel{\phi}{\rightarrow}
H_{n-g-1}\otimes S_1 \rightarrow \cdots \rightarrow H_0 \otimes
S_{n-g} \rightarrow \Gr_{n-g} \rightarrow 0\;,
\]we show that the cokernel of $\phi$, denoted by $L$, is a maximal
Cohen--Macaulay module. Because the ($n-g$)th Koszul homology
module $H_{n-g}$ is the canonical module of $S$, the exact
sequence $ 0 \rightarrow H_{n-g} \stackrel{\phi}{\rightarrow}
H_{n-g-1}\otimes S_1 \rightarrow L \rightarrow 0$ splits. On the
other hand, we have the following diagram of equivalences.

\[
\begin{CD}
\Hom_S(H_{n-g} \bigoplus L,H_{n-g}) @= \;\; \Hom_S(H_{n-g-1}
\otimes
S_1 , H_{n-g})\\
@| @|\\
\Hom_S(H_{n-g},H_{n-g}) \bigoplus \Hom_S(L, H_{n-g}) @= \bigoplus^n \Hom_S(H_{n-g-1}, H_{n-g})\\
@| @| \\
S \bigoplus \Hom_S(L,H_{n-g}) @= \bigoplus^n H_1
\end{CD}
\]

\noindent This contradicts that $H_1$ cannot have $S$ as a direct
summand (\cite[Theorem 1.4.9]{G}, \cite[Theorem 5.2.1]{V1}).
 \QED

\bigskip

Suppose that $\Gr$ is not integrally closed. Under the same
assumptions as those in Theorem~\ref{hong3}, there are examples
which show that $\Gr_{d-g-1}$ is not necessarily the first
component which is not integrally closed. On the other hand, if
$\Gr_i$ is integrally closed for all $i<d-g-1$, then we may have a
better understanding of the difference between $\Gr_{d-g-1}$ and
its bidual $\Gr_{d-g-1}^{**}=\Hom_S(\Hom_S(\Gr_{d-g-1},S),S)$.

\begin{proposition}\label{nu2}
Let $R$ be a Gorenstein local ring of dimension $d \geq 3$,
$\mathfrak{p}$ a perfect prime ideal of height $2$, generated by a
$d$--sequence. Let $\varphi$ be the matrix of syzygies of
$\mathfrak{p}$. Suppose that $S=R/\mathfrak{p}$ is an integrally
closed domain and that for every proper prime ideal $\mathfrak{q}
\supset \mathfrak{p}$,
\[\nu(\mathfrak{p}_{\mathfrak{q}})\leq \height(\mathfrak{q})-1,\]
where $\nu(\mathfrak{p}_{\mathfrak{q}})$ is the minimal number of
generators of $\mathfrak{p}_{\mathfrak{q}}$. Then
\[ \Gr_{n-2}^{**}/\Gr_{n-2} \simeq \Ext_R^d(R/I_1(\varphi), R),
\]where $\Gr_{n-2}=\mathfrak{p}^{n-2}/\mathfrak{p}^{n-1}$ and $\Gr_{n-2}^{**}=\Hom_S(\Hom_S(\Gr_{n-2},S),S)$.
\end{proposition}

\demo By Proposition~\ref{appcx1}, the associated graded ring
$\Gr=\bigoplus_{t \geq 0} \Gr_t$ of the prime ideal $\mathfrak{p}$
is a domain. We denote $\Gr_{n-2}^{**}/\Gr_{n-2}$ by $C$. From the
following two exact sequences

\[ 0 \rightarrow H_{n-2} \stackrel{\phi}{\rightarrow} H_{n-3}\otimes S_1
\rightarrow \cdots \rightarrow H_1 \otimes S_{n-3} \rightarrow
S_{n-2} \rightarrow \Gr_{n-2} \rightarrow 0, \]

\[0 \rightarrow \Gr_{n-2} \rightarrow \Gr_{n-2}^{**} \rightarrow C
\rightarrow 0,\] we obtain that

\[\Ext^{d-2}_S(C,H_{n-2}) \simeq \coker \xi,\] where $\xi:\Hom_S(H_{n-3}\otimes S_1 , H_{n-2})
\rightarrow \Hom_S(H_{n-2},H_{n-2})$ is the dual map of $\phi$. We
claim that
\[ \coker \xi \simeq R/I_1(\varphi)\;.\]

\noindent Set the matrix of syzygies $\varphi=[a_{ij}]=[v_1
\;\ldots \;v_{n-1}]$, where $1 \leq i \leq n$, $1 \leq j \leq n-1$
and $v_s$'s are the column vectors. For each $i$, we have
\[ \sum_{j=1}^{n-1} (-1)^{j-1} a_{ij} \bigwedge_{s=1\;,\;s \neq j}^{n-1} v_s  \; \in \; \mathfrak{p}\bigwedge^{n-2} R^n.\]
For each $j=1, \ldots, n-1$, the map $\phi$ is defined in the
following manner.
\[
\begin{array}{llc}
 \phi: \overline{\bigwedge_{s=1\;,\;s \neq j}^{n-1} v_s} &\mapsto
 &\sum_{t=1\;,\;t \neq j}^{n-1}(-1)^t \overline{\bigwedge_{s=1\;,\;s \neq
 j,t}^{n-1} v_s } \otimes v_t \\ & & || \\ & &
\sum_{i=1}^{n}\left( \sum_{t=1\;,\;t \neq j}^{n-1}(-1)^t
\overline{a_{it} \bigwedge_{s=1\;,\;s \neq j,t}^{n-1} v_s}
\right)e_i.
\end{array}
\] For any $h=\sum_{i=1}^{n} h_i \epsilon_i \in
\Hom(H_{n-3}\otimes S_1 , H_{n-2})$, the map $\xi$ is defined as
the following.
\[ \xi(h)=h \circ \phi: \overline{\bigwedge_{s=1\;,\;s \neq j}^{n-1} v_s} \mapsto \sum_{i=1}^{n}\left( \sum_{t=1\;,\;t \neq j}^{n-1}(-1)^t \overline{a_{it}
\bigwedge_{s=1\;,\;s \neq j,t}^{n-1} v_s \bigwedge h_i} \right).
\]Let $h_{ij}$ be $v_j \epsilon_i$, where $1
\leq i \leq n$ and $1 \leq j \leq n-1$. For example, in case when
$j=1$,
\[
\xi(h_{i1})\left(\overline{\bigwedge_{s=2}^{n-1} v_s} \right)=
\sum_{t=2}^{n-1}(-1)^t \overline{a_{it} \bigwedge_{s=1\;,\;s \neq
t}^{n-1} v_s} = a_{i1} \overline{\bigwedge_{s=2}^{n-1} v_s}.
\]Similarly we show that for each $1
\leq i \leq n$ and $1 \leq j \leq n-1$, the map $\xi(h_{ij})$ is
the multiplication by $a_{ij}$, which proves the claim. Let
$\underline{a}:a_1 , a_2$ be an $R$-sequence in $\mathfrak{p}$ and
$\underline{b}:b_1 , \ldots, b_{d-2}$ in $R$ such that
$\underline{a},\underline{b}$ is a system of parameters of $R$.
Let $J'$ be the $S$-ideal generated by the images of $b_i$'s in
$S$ for all $i=1, \ldots, d-2$. Finally we obtain the following
natural isomorphisms:

\bigskip

\[\begin{array}{lll}
C &\simeq &\Ext^{d-2}_S(\Ext^{d-2}_S(C,H_{n-2}),H_{n-2}) \\
&\simeq &\Ext^{d-2}_S(\coker \xi, H_{n-2})\\
&\simeq &\Hom_{S/J'}(\coker \xi, \omega_S/J'\omega_S)\\
&\simeq &\Hom_{S/J'}(\coker \xi, \omega_{S/J'})\\
&\simeq &\Hom_{R/(\underline{a},\underline{b})}(\coker \xi, R/(\underline{a},\underline{b}) )\\
&\simeq &\Ext^d_R(R/I_1(\varphi), R).\end{array}
\]

\QED

\begin{corollary}\label{nu1}
Let $R$ be a Gorenstein local ring, $\mathfrak{p}$ a perfect prime
ideal generated by a strongly Cohen--Macaulay $d$--sequence and
$E$ the conormal module $\mathfrak{p}/\mathfrak{p}^2$. Suppose
that $S=R/\mathfrak{p}$ is an integrally closed domain of
dimension $2$, that $\nu(\mathfrak{p})$ is $3$, and that
$\height(\mathfrak{p})$ is $2$. Then $\nu(\overline{E})$ equals to
$4$.
\end{corollary}

\begin{example}\label{ex2}{\rm
Let $R=k[x,y,z,u,v]$ and let $A$ be the matrix of syzygies of
$R$--ideal $\mathfrak{p}$:
\[
A=\left[\begin{array}{lll} x&v-u&z\\y&x&v\\z&u&x\\u&z&y
\end{array}\right]
\]

\noindent The associated graded ring $\Gr=\gr_{\mathfrak{p}}
R=\bigoplus \Gr_i$ is a domain but it is not integrally closed.
The conormal module $E=\Gr_1=\mathfrak{p}/\mathfrak{p}^2$ is
reflexive. But $\nu(\Gr_2)$ and $\nu(\Gr_2^{**})$ are $10$ and
$11$ respectively. This shows that $\Gr_2$ is not integrally
closed because of the failure of the $\mathfrak{m}$-fullness. In
particular, by Proposition~\ref{nu2}, $\Gr_2^{**}$ is the integral
closure of $\Gr_2$. }\end{example}

\end{document}